\documentclass[amscd,amssymb,verbatim,12pt]{amsart}
\usepackage{epsfig}
\usepackage{graphicx}

\usepackage{pdfpages}

\newif\ifAMS
\IfFileExists{amssymb.sty}
  {\AMStrue\usepackage{amssymb}}
  {\usepackage{latexsym}}
  \usepackage{enumerate}
\theoremstyle{plain}

\newtheorem{Thm}{Theorem}[section]
\newtheorem{Cor}[Thm]{Corollary}
\newtheorem{Lem}[Thm]{Lemma}
\newtheorem{Prop}[Thm]{Proposition}

\theoremstyle{definition}
\newtheorem{Def}{Definition}

\theoremstyle{remark}

\newtheorem{Qu}[Thm]{Question}

\DeclareMathOperator{\length}{length}
\DeclareMathOperator{\diam}{diam}
\DeclareMathOperator{\asdim}{asdim}

\newcommand{\interior}{^{ \kern-5pt ^\circ}}

\newcommand {\al}{\alpha}
\newcommand {\be}{\beta}

\newcommand {\R}{{\mathbb R}}

\begin{document}
\title
{Asymptotic dimension of planes and planar graphs}

\author
{Koji Fujiwara }
\email{kfujiwara@math.kyoto-u.ac.jp}
\address{Department of Mathematics, Kyoto University,
Kyoto, 606-8502, Japan}

\author
{Panos Papasoglu }


\email {} \email {papazoglou@maths.ox.ac.uk}

\address
{Mathematical Institute, University of Oxford, 24-29 St Giles',
Oxford, OX1 3LB, U.K.  }

\address
{ }

 \thanks{ The first author is
    supported in part by Grant-in-Aid for Scientific Research
    (No. 15H05739, 20H00114).
    He would like to thank Mathematical Institute of University of 
    Oxford for their hospitality.}

\begin{abstract} 
We show that the asymptotic dimension of a geodesic space that 
is homeomorphic to a subset in the plane
is at most three. In particular, the asymptotic dimension of the plane
and any planar graph is at most three.

\end{abstract}
\maketitle

\section{Introduction}
\subsection{Statements}
The notion of {\it asymptotic dimension} introduced by Gromov \cite{Gr} has become central in Geometric Group Theory mainly because of its
relationship with the Novikov conjecture. 
The asymptotic dimension ${\rm asdim}\,X$ of a metric space $X$ is defined as follows: ${\rm asdim}\,X \leq n$ if and only if for every $m > 0$ there exists $D(m)>0$ and a covering $\mathcal{U}$ of $X$ by sets of diameter $\leq D(m)$ 
($D(m)$-bounded sets) such that any $m$-ball in $X$ intersects at most $n+1$ elements of $\mathcal{U}$.
We say ${\rm asdim} X \le n$, {\it uniformly} if one
can take $D(m)$ independently from $X$ if it belongs to a certain family.

In this paper we deal with asymptotic dimension in a purely geometric setting, that of Riemannian planes and planar graphs.
An aspect of the geometry of Riemannian planes that is studied extensively is that of the isoperimetric problem-even though in that case one
usually imposes some curvature conditions (see \cite{BC},\cite{MHH},\cite{HHM}, \cite{R}, \cite{H}, \cite{GP}).
We note that Bavard-Pansu (\cite{BP}, see also \cite{Bo})
have calculated the minimal volume of a Riemannian plane. There are some general results in the related case of a 2-sphere \cite{He}. On the other hand there is a vast literature dealing with planar graphs.
See eg \cite{AH},\cite{GHT},\cite{LT},\cite{NC},\cite{Tu}.

We prove the following:

\begin{Thm}\label{Main}
Let $P$ be a geodesic metric space that is homeomorphic to $\Bbb R^2$. Then the asymptotic dimension of $P$ is at most three, uniformly.
More generally if $P$ is a geodesic metric space such that
there is an injective continuous map from $P$ to $\Bbb R^2$, 
then the conclusion holds. 

To be more precise,  the following holds: Given $m>0$ there is some $D(m)>0$ such that  there is a cover of $P$ with sets of diameter $<D(m)$ and that
any ball of radius $m$ intersects at most 4 of these sets.

Moreover, we can take $D(m)=3 \cdot 10^6m$.

\end{Thm}
We note that any complete Riemannian metric
on $\R^2$ gives an example of such a geodesic space $P$.

We say a connected graph $P$ is  {\it planar} 
if there is an injective map
$$\phi:P \to \Bbb R^2$$
such that on each edge of $P$, the map $\phi$ is continuous.

We view a connected graph as a geodesic space
where each edge has length $1$. We denote this metric by $d_P$.
We do not assume that the above map $\phi$ is continuous on $P$
with respect to $d_P$, so that Theorem \ref{Main}
might not directly apply, but 
the same conclusion holds for planar graphs.

\begin{Thm} \label{inf-graph}
The asymptotic dimension of a planar graph, $(P,d_P)$, is at most three, uniformly
for all planar graphs. 
\end{Thm}
The conclusion on the existence of a covering in Theorem \ref{Main}
holds for planar graphs as well. 

The proof of both theorems will be given in Section \ref{section.asdim.three}.

There is a notion called {\it Assouad-Nagata dimension}, which 
is closely related to asymptotic dimension. The only difference is that 
it additionally requires that there exists a constant $C$ such that 
$D(m) \le C m$ in the definition of asymptotic dimension.
Since we have a such bound, we also prove that 
Assouad-Nagata dimension of $P$ is at most three
in Theorems \ref{Main} and \ref{inf-graph}.

We note that all finite graphs have asymptotic dimension 0 however our theorem makes sense
for finite graphs as well. We restate Theorem \ref{inf-graph}
in terms of a covering 
for finite planar graphs as a special case:

\begin{Cor}\label{fin-graph}
For any $m>0$ there is $D(m)>0$ such that if $G$ is any finite planar graph there is a cover of $G$
by subgraphs $G_i, i=1,...,n$ such that the diameter of each $G_i$ is bounded by $D(m)$ and any ball
of radius $m$ intersects at most 4 of the $G_i$'s.
\end{Cor}

In connection to Theorem \ref{inf-graph}, we would like to mention the following theorem.
\begin{Thm}[Ostrovskii-Rosenthal]\cite{OR}
If $\Gamma$ is a connected graph with finite degrees 
excluding the complete graph $K_m$ as a minor, then $\Gamma$
has asymptotic dimension at most $4^m-1$.
\end{Thm}
$K_m$ here is the compete graph of $m$-vertices.
The degree of a vertex is the number of edges incident at the vertex.
A {\it minor} of a graph $\Gamma$ is a graph $M$ obtained by contracting edges in a subgraph of $\Gamma$. 
The well-known Kuratowski Theorem states that
a finite graph is planar if and only if the $K_5$
and $K_{3,3}$, the complete bipartite graph on six vertices,
are excluded as {\it  minors} of the graph.
This characterization applies to infinite graphs
if one defines an infinite graph to be planar provided
there is an embedding of the graph into $\Bbb R^2$, \cite{DS}.
So, as a special case, the theorem above implies that
an infinite finite degree graph  that embeds
in $\Bbb R^2$ has asymptotic dimension 
at most $4^5-1$, in particular finite. 
We also remark that they also proved this bound 
for Assouad-Nagata dimension, which bounds 
asymptotic dimension from above. The proof relies on earlier results of
Klein, Plotkin, and Rao \cite{KPR}.

\subsection{Idea of proofs}
We give an outline of the proof of our results. 
We fix a basepoint $e$ in $P$ and we consider `annuli' around $e$ of a fixed width
(these are metric annuli so, if $P$ is a plane with a Riemannian metric, topologically
are generally discs with finitely many holes). 
Here, annuli are subsets defined as follows:
Consider $f(x)=d(e,x)$.
Fix $m>0$. We will pick $N\gg m$ and consider for $k\in \mathbb N$
the ``annulus''
$$A_k(N)=\{x|kN \le f(x) < (k+1)N\}$$

We show in section 3 that in the large scale these annuli resemble cacti. Generalizing a well known result
for trees and $\mathbb R$-trees we show in section 2 that cacti have asymptotic dimension at most 1. We show in section 3 that `coarse cacti' also have asymptotic
dimension 1. In section 4 we decompose our space in `layers' which are coarse cacti which implies that the asymptotic dimension of the space is at most 3.

In the proofs in sections 2-4 the constants and inequalities that we use are far from optimal, we hope instead that they are `obvious' and easily verifiable by the reader.

In section 5 we show that our result can not be extended to Riemannian metrics on $\mathbb R^3$ and we pose some questions.
We give some updates as notes added in proof.

 \section*{Acknowledgements}
 We thank Romain Tessera for his comments and Agelos Georgakopoulos for bringing \cite{OR} to our attention.
 We thank Urs Lang for letting us know the work \cite {JL}.
 We are grateful to the referee for very carefully reading the manuscripts
 and making precise and insightful comments.

\section{Asymptotic dimension of cacti}

\subsection{Cactus}
As we said, the idea of our proof is that the successive `annuli' making up the plane resemble cacti and so they have
asymptotic dimension at most 1. 

We begin by showing that a cactus has asymptotic dimension at most 1.

\begin{Def}[Cactus]
A \textit{cactus} (graph) is a connected graph such that any two cycles intersect at at most one point.
More generally we will call cactus a geodesic metric space $C$ such that any two distinct simple closed
curves in $C$ intersect at at most one point.
\end{Def}

We remark that our  notion of cactus generalizes the classical graph theoretic notion
in a similar way as $\mathbb R$-trees generalize trees.
Historically, a cactus graph was introduced by K. Husimi and studied in \cite{HU}.
Cacti have been studied and used in graph theory, algorithms,  electrical engineering and others.

\begin{Prop}\label{cactus.asdim1}
A cactus $C$ has ${\rm asdim}\, \le 1$,  uniformly over all cacti. 
Moreover, we can take $D(m)=1000m$.

\end{Prop}

\proof
Let $m>0$ be given. It is enough to show that there is a covering of $C$ by uniformly bounded sets
such that any ball of radius $m$ intersects at most 2 such sets.
Fix $e \in C$. Consider $f(x)=d(e,x)$.
We will pick $N=100 m$ and consider for $k\in \mathbb N \cup \{0\}$
the ``annulus''
$$A_k=\{x|kN \le f(x) < (k+1)N\}.$$
We define an equivalence relation on $A_k$: $x\sim y$ if there are $x_1=x,x_2,...,x_n=y$ such that $x_i\in A_k$
and $d(x_i,x_{i+1})\leq 10m$ for all $i$. Since every $x\in C$ lies in exactly one $A_k$ this equivalence relation
is defined on all $C$. Let's denote by $B_i$, $(i\in I)$ the equivalence classes of $\sim $ for all $k$. By definition, for each $A_k$,  if $B_i,B_j$ lie
in  $A_k$ then a ball $B$ of radius $m$ intersects at most one of them. It follows that a ball of radius $m$
can intersect at most two equivalence classes. So it suffices to show that the $B_i$'s are uniformly bounded.
We claim that $\diam (B_i)\leq 10N$. This will show we can take
$$D(m)=1000m.$$
 We will argue by contradiction: let $x,y\in B_i\subseteq A_k$ such that $d(x,y)>10N$. We will show 
that there are two non-trivial loops on $C$ that intersect along a non-trivial arc.

Let $\gamma _1,\gamma _2$ be geodesics from
$e$ to $x,y$ respectively. Let $p$ be the last intersection point of $\gamma _1, \gamma _2$. 
We may assume without loss of generality that $\gamma _1\cap \gamma_ 2$ is an arc
with endpoints $e,p$.


By the definition of $\sim $ there is a path $\alpha $ from $x$ to $y$ that lies in the $10m$-neighborhood of $A_k$.
We may assume that $\alpha $ is a simple arc and that its intersection with each one of $\gamma _1, \gamma _2$ is connected.
If $x_1$ is the last point of intersection of $\alpha $ with $\gamma _1$ and $y_1$ is the first point of intersection of $\alpha $ with $\gamma _2$
then the subarcs of $\gamma _1, \alpha , \gamma _2$ with endpoints respectively $p,x_1$, $x_1,y_1$, $y_1,p$ define a simple closed curve $\beta $.
We note that $$d(e,x_1)\geq \length (\gamma _1)-N-10m, d(e,x_2)\geq \length (\gamma _2)-N-10m .$$
Let $\alpha _1$ be the subarc of $\alpha $ with endpoints $x_1,y_1$.  Then $$\length (\alpha _1)\geq 7N.$$ Let $x_2$ be the midpoint of $\alpha _1$.
%
\begin{figure}[htbp]
\hspace*{-3.3cm}     
\begin{center}
                                                      
   \includegraphics[scale=0.500]{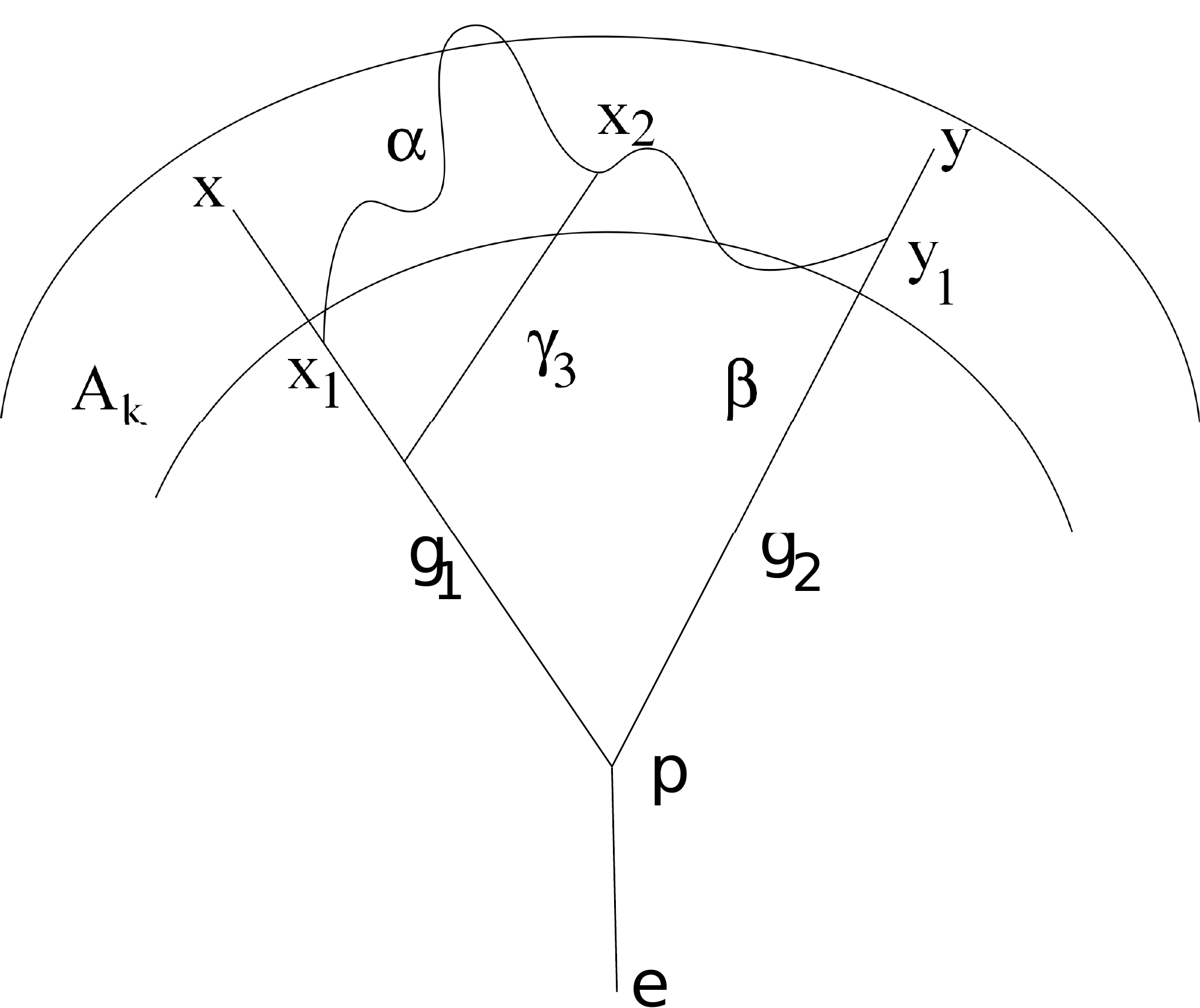}%
   \end{center}

\caption{Two loops intersecting along an arc}
  \label{}
\end{figure}

We consider a geodesic $\gamma _3$ joining $e$ to the midpoint $x_2$ of $\alpha _1$.
 We may and do assume $\gamma_1 \cap
\gamma_2 \cap \gamma_3$ is connected.
We note that $\gamma _3$ is not contained in $\beta \cup (\gamma _1\cap \gamma _2)$. Indeed if it were contained in this union then we would have,
for at least one of $i=1,2$, 
$$\length (\gamma _3)\geq \length (\gamma _i)+2N \text{ for }i=1 \text{ or } 2,$$ however this is impossible since for both $i=1,2$ we have$$d(e,x_2)\leq \length (\gamma _i)+N+10m.$$

Therefore there are two cases:

\textit{Case 1}. There is a subarc of $\gamma _3$ with one endpoint $a_1$ on $\gamma _1\cap \gamma _2$ and another endpoint $a_2\ne p$ on $\beta $
which intersects $\gamma _1\cup \beta $ only at its endpoints.
In this case we consider the loop $\gamma $ consisting of the arc on $\gamma _3$ with endpoints $a_1,a_2$ and a simple arc on $\gamma _1\cup \beta$
joining $a_1,a_2$. Clearly $\gamma $ intersects $\beta $ along a non-trivial arc contradicting the fact that $C$ is a cactus.

\textit{Case 2}. There is a subarc of $\gamma _3$ with endpoints $a_1,a_2$ on $\beta $ which intersects $\beta $ only at its endpoints.
In this case we consider the loop $\gamma $ consisting of the arc on $\gamma _3$ with endpoints $a_1,a_2$ and a simple arc on $\beta$
joining $a_1,a_2$. Clearly $\gamma $ intersects $\beta $ along a non-trivial arc contradicting the fact that $C$ is a cactus.

The moreover part follows since for a given $m>0$, we chose
$N=100m$ and showed $\diam(B_i) \le 10N$, which does not depend
on the cactus $C$.
\qed

The following is immediate from Proposition \ref{cactus.asdim1}.
\begin{Cor}
 If $X$ is quasi-isometric to a cactus then ${\rm asdim}\, X \le 1$.
Moreover if $X$ is uniformly quasi-isometric to a cactus, 
then ${\rm asdim}\, X \le 1$, uniformly.

%
\end{Cor}
To be concrete, the conclusion says that $D(m)$ in the 
definition of the asymptotic dimension depends
only on $m$ and the quasi-isometry constants.

\section{Coarse cacti}
We prove now that if a space looks coarsely like a cactus it has asymptotic dimension at most 1. We make
precise what it means to look coarsely like a cactus below.

\begin{Def}[$M$-fat theta curve]
Let $X$ be a geodesic metric space. Let $\Theta $ be a unit circle in the plane together with a diameter.
We denote by $x,y$ the endpoints of the diameter and by $q_1,q_2,q_3$ the 3 arcs joining them (ie the closures of the connected components of $\Theta \setminus \{x,y\}$).
A \textit{theta-curve} in $X$ is a continuous map $f:\Theta \to X$. Let $p_i=f(q_i),\, i=1,2,3,\, a=f(x),b=f(y)$.

A theta curve is $M$-\textit{fat} if 
there are arcs $\al _i,\be _i\subseteq p_i,\, i=1,2,3$ where $a\in \al _i,b\in \be _i$
so that 
the following hold:
\begin{enumerate}
\item
If $p_i'=p_i\setminus \al _i\cup \be _i$ then $p_i'\ne \emptyset $ and for
any $i\ne j$ and any $t\in p_i',s\in p_j'$ we have $d(t,s)\geq M$.
\item
$p_i'\cap \al _j=\emptyset,\, p_i'\cap \be _j=\emptyset $ for all $i,j$ (note by definition $p_i'$ is an open arc,
ie does not contain its endpoints).
\item
For any $t\in \al _1\cup \al _2\cup \al _3, s\in  \be _1\cup \be _2\cup \be _3$, we have $d(t,s)\geq 2M$.
\end{enumerate}
We say that $a,b$ are the \textit{vertices} of the theta curve.
We say that the theta curve is \textit{embedded} if the map $f$ is injective.
We will often abuse notation and identify the theta curve with its image giving simply the arcs of the theta curve.
So we will denote the theta curve defined above by $\Theta (p_1,p_2,p_3)$.

\end{Def}

We note that if $i\not=j,k$ then 
$$p'_i \backslash N_{M}(p_j \cup p_k) \not= \emptyset,$$
where $N_a(B)$ denotes the open $a$-neighborhood of $B$.
This is immediate from the definition. 
Indeed, let $z \in p_i'$ be a point with $d(x,\alpha_1\cup\alpha_2 \cup \alpha_3 \cup
\beta_1 \cup \beta_2 \cup \beta_3) \ge M$.
Such $z$ exists by the property (3).
But then, $d(z,p_j') \ge M$ and $d(z,p_k') \ge M$ by (1), which
implies $z \in p'_i \backslash N_{M}(p_j \cup p_k) $.

We remark that to show that a theta curve $\Theta (p_1,p_2,p_3)$ is $M$-fat it is enough to specify arcs
$p_i'\subset p_i, i=1,2,3$ so that the conditions 1,2,3 of the definition above hold. In other words the arcs $p_i'$
determine the arcs $\al _i,\be _i$.

Note that theta curves are not necessarily embedded. However we have the following:
\begin{Lem}\label{embedding}
Suppose a geodesic space $(A,d_A)$ contains an $M$-fat theta
curve $\Theta(p_1,p_2,p_3)$. Then $A$ contains an embedded $M$-fat theta
curve
$\Theta(\gamma _1,\gamma _2,\gamma _3)$, which is a subset of $\Theta(p_1,p_2,p_3)$.
\end{Lem}

\proof
Let $a,b$ be the vertices of $\Theta(p_1,p_2,p_3)$ and let $\al _i,\be _i\subseteq p_i,\, i=1,2,3$ where $a\in \al _i,b\in \be _i$
arcs as in the definition of $M$-fat theta curve. We may replace each of $p_i'=p_i\setminus \al _i\cup \be _i$ by a simple arc,
with endpoints say $a_i,b_i$. Similarly we may replace each of $\al _i,\be _i$ by simple arcs with the same endpoints.

Let $c_2,c_3$ be the last points, along $\alpha_1$ from 
$a$ to $a_1$,  of intersection of $\alpha _1,\alpha _2$ and $\alpha _1,\alpha _3$ respectively.

If $\alpha$ is an arc we denote below by $\alpha (u,v)$ the subarc of $\alpha $ with endpoints $u,v$.

\begin{figure}[htbp]\label{2}
\begin{center}
\includegraphics[scale=0.45]{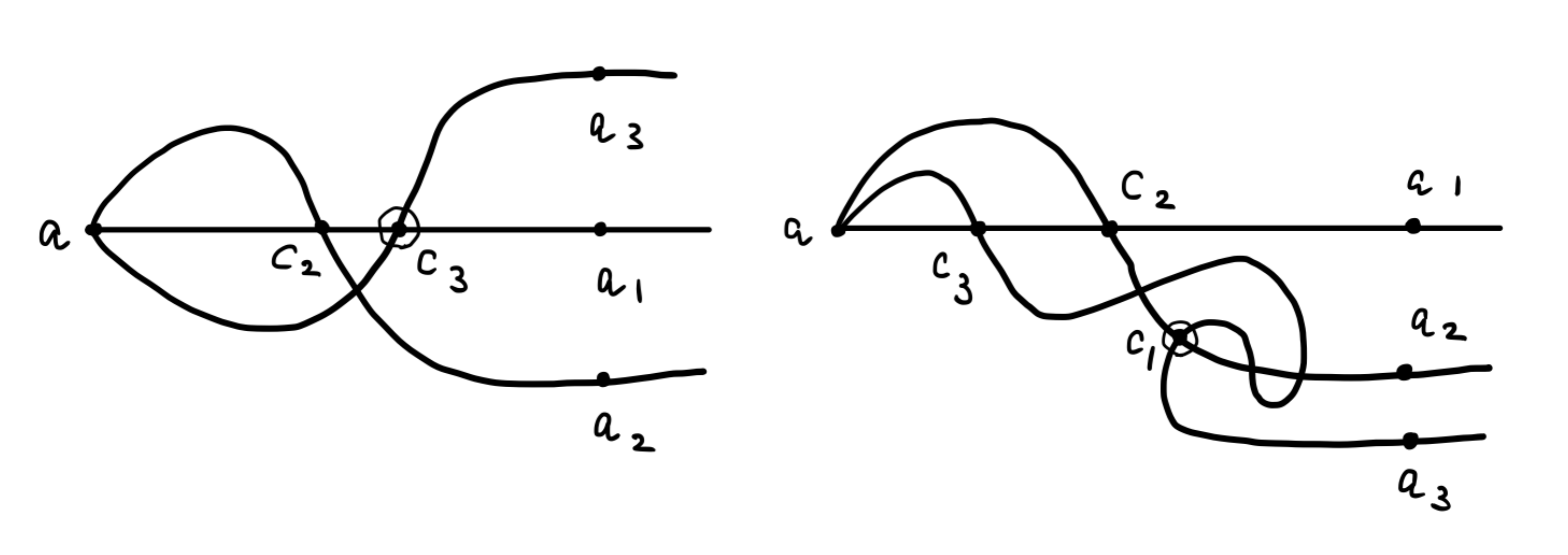}
  \end{center}
  \caption{Left for Case 1 and right for Case 2.}
  \label{fig.embedding}
\end{figure}

We divide the case into two depending on the position of 
$c_2,c_3$ on $\alpha_1$.
See Figure \ref{fig.embedding}.

(i) 
Suppose $c_3\in \alpha _1(c_2,a_1)$.
We further divide the case into two: 

{\it Case 1}. $(\alpha_3(c_3,a_3)\backslash c_3)\cap
(\alpha_2(c_2,a_2)\backslash c_2) =\emptyset$.
Then, we take $c_3$ to be a vertex of the new theta curve and replace $\al _i, i=1,2,3$ 
by $$\alpha _1(c_3,a_1), \, \, \alpha_1(c_3,c_2) \cup \alpha _2(c_2,a_2), \, \, 
\alpha _3(c_3,a_3).$$

{\it Case 2}.
$(\alpha_3(c_3,a_3)\backslash c_3)\cap
(\alpha_2(c_2,a_2)\backslash c_2) \not=\emptyset$.
Then, let $c_1$ be the last point, along $\alpha_3$, of the 
intersection $\alpha_3(c_3,a_3) \cap
\alpha_2(c_2,a_2)$.
In this case, we take $c_1$ to be a vertex of the new theta curve and 
replace  $\al _i, i=1,2,3$ by
$$\alpha_2(c_1,c_2) \cup \alpha_1(c_2,a_1), 
\alpha_2(c_1, a_2),
\alpha_3(c_1, a_3).$$

(ii) Suppose $c_3 \in \alpha_1(a,c_2)$.
In this case, we replace $\alpha_i$ with $\alpha_i'$
after we switch the roles of $\alpha_2$ and $\alpha_3$, 
so that $c_2$ and $c_3$ are switched and we are in (i).

In all cases, any pair of  $\alpha_i'$ intersect only in the new vertex, and 
$(\alpha_1' \cup \alpha_2' \cup \alpha_3')
\subset
(\alpha_1 \cup \alpha_2 \cup \alpha_3)$.


We replace $\be _i$ similarly. Clearly we obtain in this way an $M$-fat embedded theta curve.
\qed

\begin{Def}[$M$-coarse cactus]
Let $X$ be a geodesic metric space. If there is an $M>0$ such that $X$ has no embedded, $M$-fat theta curves then we say that
$X$ is an $M$-\textit{coarse cactus} or simply a \textit{coarse cactus}.
\end{Def}

We give now a proof that a coarse cactus
has asymptotic dimension at most one imitating the proof of Proposition\ref{cactus.asdim1}.

\begin{Thm}\label{quasi.cactus}
Let $C$ be an $M$-coarse cactus. 
Then $\asdim C \le 1$.
Moreover, it is uniform with $M$ fixed. 
Further, for any $m\ge M$, we can take
$D(m)=10^5m$.
\end{Thm}
Note that, for $m < M$, we could put, for example,  $D(m)=10^5M$, so that 
we can set $D(m)=10^5 \max\{m,M\}$ for all $m$.
\proof

Let $m>0$ be given.  
It is enough to show that there is a covering of $C$ by uniformly bounded sets
such that any ball of radius $m$ intersects at most 2 such sets. 
Without loss of generality we may assume $m\ge M$.
Fix $e \in C$. Consider $f(x)=d(e,x)$.
We will pick $N=100m$ and consider 
the ``annulus''
$$A_k=\{x|kN \le f(x) < (k+1)N\}.$$
We define an equivalence relation on $A_k$: $x\sim y$ if there are $x_1=x,x_2,...,x_n=y$ such that $x_i\in A_k$
and $d(x_i,x_{i+1})\leq 10m$ for all $i$. Since every $x\in C$ lies in exactly one $A_k$ this equivalence relation
is defined on all $C$. Let's denote by $B_i$, $(i\in I)$ the equivalence classes of $\sim $. By definition if $B_i,B_j$ lie
in some $A_k$ then a ball $B$ of radius $m$ intersects at most one of them. It follows that a ball of radius $m$
can intersect at most two equivalence classes. So it suffices to show that the $B_i$'s are uniformly bounded.
We claim that $\diam (B_i)\leq 1000N$, which shows it suffices to take 
$$D(m)=1000N=100000m.$$
We will argue by contradiction: let $x,y\in B_i\subseteq A_k$ such that $d(x,y)>1000N$. We will show 
that there is an $N$-fat theta curve in $C$,
which is a contradiction since $N >M$, and Lemma \ref{embedding} applies.

Since $\diam A_k \le 2(k+1)N$, we may assume $k\ge 499$, so that 
$d(e,x) \ge 499N$ for $x\in A_k$.

Let $\gamma _1:[0,\ell_1]\to C,\gamma _3:[0,\ell_3]\to C$ be geodesics (parametrized with respect to arc length) from
$e$ to $x,y$ respectively.

By the definition of $\sim $ there is a path $\alpha:[0,\ell]\to C $ from $x$ to $y$ that lies in the $10m$-neighborhood of $A_k$. 
We further assume that $\alpha $ is simple.
 Let $a \in  \alpha$ such that $$d(a,x)=d(a,y).$$
 Note $d(a, x)=d(a,y) > 500 N$.

We consider a geodesic $\gamma _2:[0,\ell_2]\to C$ joining $e$ to $a$. 
We claim that the theta curve $$\Theta=\gamma_1\cup \gamma _2\cup \gamma _3 \cup \alpha $$
with vertices $e,a$
is $N$-fat.  Explicitly the 3 arcs of $\Theta $ are $p_1=\gamma _1\cup \alpha (x,a)$, $p_2=\gamma _2$
and $p_3=\gamma _3 \cup \alpha (a,y)$.

To see that $\Theta $ is $N$-fat it is enough to define subarcs $p_i'\subseteq p_i$ so that the conditions of
the definition of $N$-fat theta curves are satisfied. 

We set $p_i'=\gamma _i[\ell_i-20N,\ell_i-10N], i=1,2,3$. We follow the notation of the definition of $M$-fat theta curve, and we denote by
$\alpha _i, \beta _i$ ($i=1,2,3)$ the arcs of the theta curve containing $a,e$ respectively.
We verify the properties (1), (2), (3).
Note that $(k-1)N \le \ell_i \le (k+2)N$, and $499 \le k$.
Also,
$$\alpha \cap (p_1' \cup p_2' \cup p_3' \cup \beta_1 \cup \beta_2\cup \beta_3)=\emptyset.$$

\begin{figure}[htbp]\label{2}
\begin{center}
\includegraphics[scale=0.4]{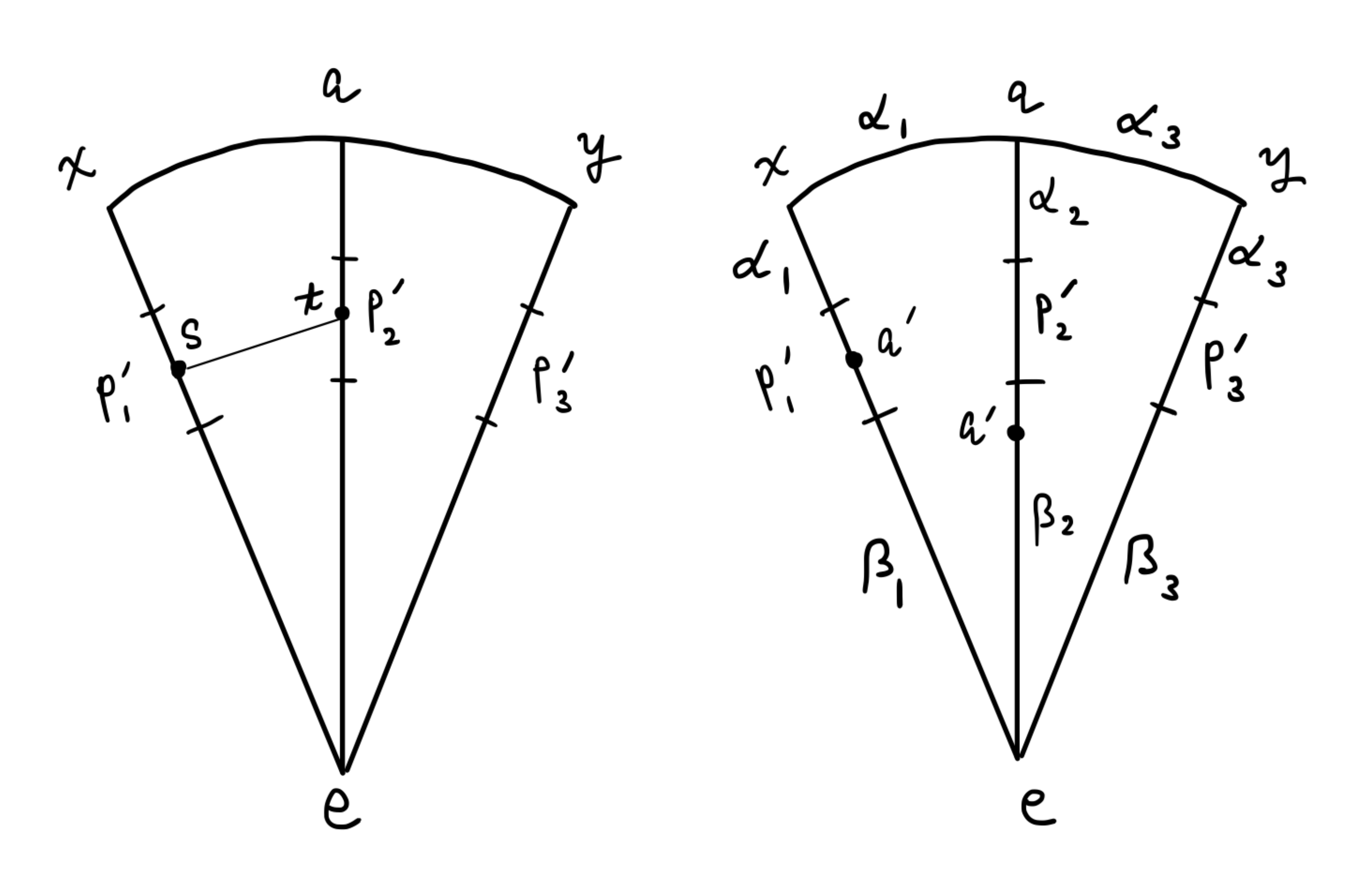}
  \end{center}
  \caption{Left figure for (1) and  right figure for (2)}
  \label{nofattheta}
\end{figure}

(1).If there are $s\in p_i',t\in p_j'$ such that $d(s,t)<N$ then 
it follows, by the triangle inequality,  that $d(a,x)<20N+N+20N=41N$ or $d(a,y)<41N$ or $d(x,y)<41N$ which is a contradiction 
since
$d(a,x)>500N$, $d(a,y)>500N$, and $d(x,y)>1000N$.
See Figure \ref{nofattheta}.

(2). In the case of $p_1'$, $p_1' \cap \alpha_1 =\emptyset, p_1' \cap \beta_1=\emptyset$ is 
trivial by definition. 
If $p_1' \cap \alpha_2\not=\emptyset$,  then
$d(x,a) \le 20N+10N=30N$, impossible.
If $p_1' \cap \alpha_3 \not=\emptyset$ 
then, $p_1' \cap \alpha = \emptyset$ implies that 
$p_1' \cap (\gamma_3 \cap \alpha_3) \not=\emptyset$, 
so that 
$d(x,y) \le 20N + 10N=30N$, impossible.
If $p_1' \cap \beta_2 \not=\emptyset$, 
then let $a' \in \gamma_2$ be a point in the intersection. See Figure \ref{nofattheta}.
Then $d(a,a') \le 22N$.
This is because
$$d(e,x) = d(e,a')+d(a',x) \le d(e,a) -d(a,a') + 20N, $$
but since $|d(e,x)-d(e,a)|\le 2N$, we conclude $d(a,a') \le 22N$.
Therefore $d(a,x) \le d(a,a') +d(a',x) \le 22N + 20N =42N$, impossible. 
If $p_1' \cap \beta_3 \not=\emptyset$, then $d(x,y) \le 42N$, impossible.
We are done with $p_1'$.

In the case of $p_2'$.
$p_2' \cap \alpha_2 =\emptyset, p_2' \cap \beta_2=\emptyset$ is 
trivial. 
If $p_2' \cap \alpha_1 \not=\emptyset$, then 
$d(a,x) \le 20N+10N=30N$, impossible (use $p_2' \cap \alpha=\emptyset$).
Same for $p_2' \cap \alpha_3=\emptyset$.
If $p_2' \cap \beta_1 \not=\emptyset$, then 
as we argued for $p_1' \cap \beta_2 =\emptyset$,
we would have $d(a,x) \le 42N$, impossible.
The argument is same for $p_2' \cap \beta_3 =\emptyset$.
Therefore the condition holds for $p_2'$.

In the case of $p_3'$.
The argument is exactly same  as $p_1'$.

(3). 
If $t\in \alpha$, then $d(e,t) \ge (k-1)N$. If $t\in \alpha_i \cap \gamma_i$ for some $i$, then 
$$d(e,t) \ge \ell_i -10N \ge kN-11N.$$
So, if $t\in \alpha_1\cup \alpha_2 \cup \alpha_3$, then $d(e,t) \ge kN-11N$.
On the other hand, if $s\in \beta_i$ for some $i$, then 
$$d(e,s) \le \ell_i -20N \le kN-18N.$$
It follows that $d(t,s) \ge 7N \ge 700M$.
This completes the proof. 
\qed

We conclude this section with a lemma that is a consequence 
of the Jordan-Schoenflies curve theorem.

\begin{Lem}[The theta-curve lemma]\label{theta-curve-lemma}
Let $\Theta(p,q,r)$ be an embedded theta curve in $\Bbb R^2$, and 
$e\in \Bbb R^2$  a point with $e \not\in \Theta$.
Then after swapping the labels $p,q,r$ if necessary, the simple loop $p \cup r$
divides $\Bbb R^2$ into two regions  such that one contains $e$ and
the other contains (the interior of) $q$.

\end{Lem}

\proof
By the Jordan-Schoenflies curve theorem (cf. \cite{C}), after applying a self-homeomorphism
of $\Bbb R^2$, we may assume the simple loop $p\cup r$ is the unit circle in $\Bbb R^2$, which divides the plane into two regions, $D_1, D_2$. 
If $e$ and $q$ are not in the same region, we are done. 
So, suppose  both are in, say, $D_1$. Then the arc $q$ divides
$D_1$ into two regions, and call the one that contains $e$, $D_1'$. 
After swapping $p,r$ if necessary, the boundary of $D_1'$ 
is the simple loop $p \cup q$.
Now, apply the Jordan-Schoenflies curve theorem to the loop
$p \cup q$, then it divides the plane into two regions such that 
one is $D_1'$ and the other one contains $r$.
Finally we swap $q,r$ and we are done. 
\qed

\section{Asymptotic dimension of planar sets and graphs}\label{section.asdim.three}

\begin{Def}[Planar sets and graphs]
Let $(P,d_P)$ be a geodesic metric space.
We say it is a {\it planar set} if 
there is an injective continuous map,
$$\phi: P \to \Bbb R^2.$$

Let $P$ be a graph. We say $P$ is  {\it planar} 
if there is an injective map
$$\phi:P \to \Bbb R^2$$
such that on each edge of $P$, the map $\phi$ is continuous.
\end{Def}
We view a connected graph as a geodesic space
where each edge has length $1$. We denote this metric by $d_P$.
We do not assume that the above map $\phi$ is continuous 
with respect to $d_P$ when $P$ is a graph. 

\subsection{Annuli are coarse-cacti}
Let $(P,d_P)$ be a geodesic metric space and pick a base point $e$.
For $r>m>0$, set
$$A(r,r+m)=\{x \in P| r \le d_P(e,x) < r+m\},$$
which we call an {\it annulus}, although it is not always
a topological annulus.

We start with a key lemma. 


\begin{Lem}\label{annulus.no.theta}
Suppose $(P,d_P)$ is a planar set or a planar graph. 
Then, for any $r,m>0$, each connected component, $A$, of 
$A(r,r+m)$ with the path metric $d_A$ has no 
embedded $m$-fat theta curve. 
\end{Lem}

\proof
{\it Case 1: Planar sets}.
We argue by contradiction.
Suppose  $A$ contains an embedded  $m$-fat theta-curve
$\Theta(p,q,s)$.

As we noted after the definition of a fat theta curve (recall $p' \subset p$):
$$p \setminus N_m(q\cup s) \not=\emptyset,
q \setminus N_m(s \cup p) \not=\emptyset,
s \setminus N_m(p\cup q) \not=\emptyset.
$$
Here, $N_m$ is for the open $m$-neighborhood w.r.t. $d_A$.

Using the map $\phi$, we can identify $P$ with its image in $\Bbb R^2$.
Since $\Theta$ is (continuously) embedded by $\phi$, we view it as a subset
in $\Bbb R^2$.
Then by the theta-curve lemma (Lemma \ref{theta-curve-lemma}), after swapping $p,q,s$ if necessary, 
the simple loop $p \cup s$ divides $\Bbb R^2$ into two 
regions such that one contains $e$ and the other contains (the interior of)
the arc $q$.

Take a point 
$$x \in q\setminus N_m(s\cup p).$$

\begin{figure}[htbp]
\begin{center}
\includegraphics[scale=0.2, angle=-90]{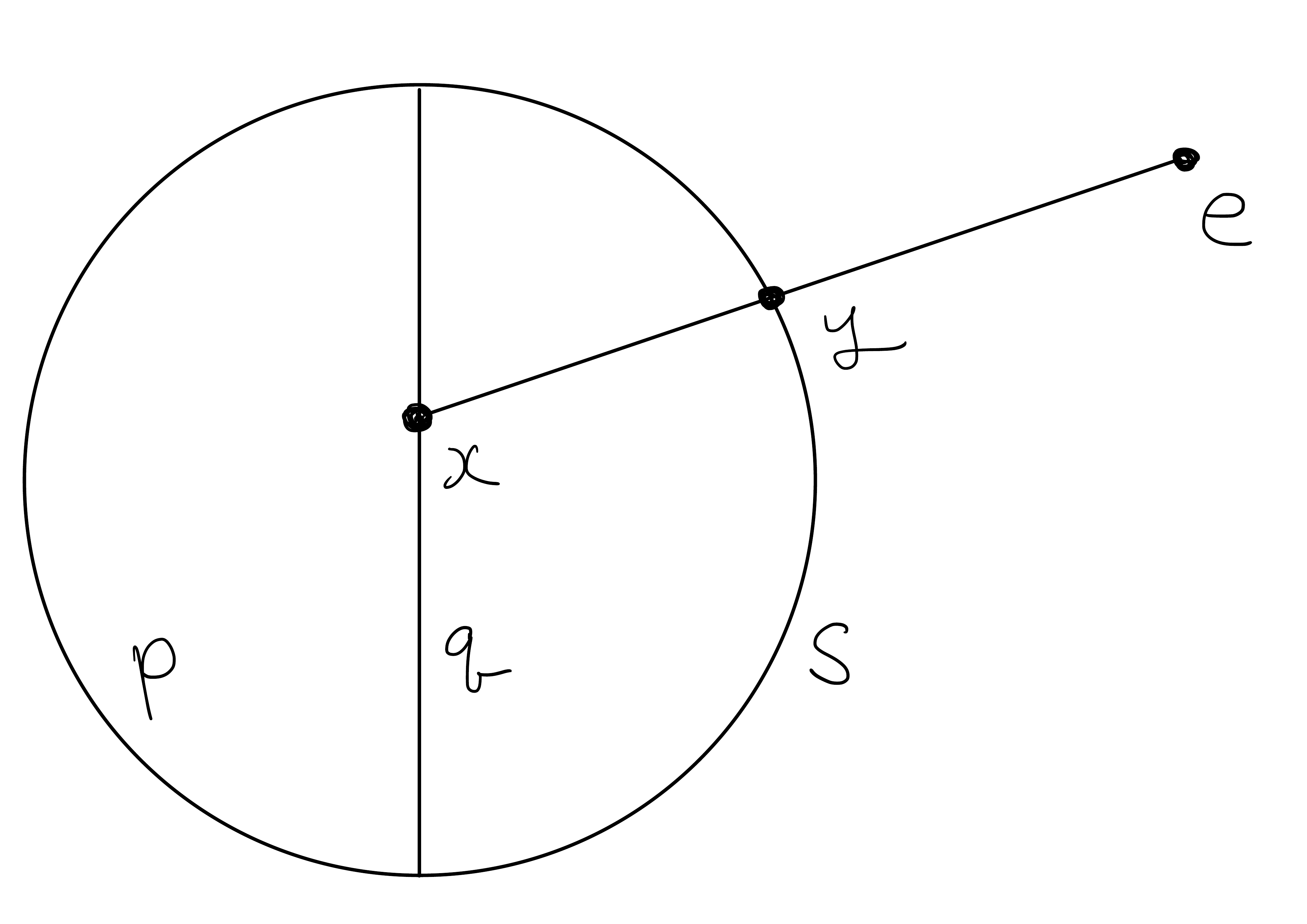}
  \end{center}
  \caption{$\gamma=[e,x]$ must intersect $p \cup s$}
  \label{notheta}
\end{figure}

Join $e$ and $x$ by a geodesic $\gamma$ in 
the space $P$. 
Then by the Jordan curve theorem, $\gamma $ must intersect 
$p \cup s$ since $x \not\in D$.  
See Figure \ref{notheta}.

Let $y$ be a point on $\gamma$ that is on $p \cup s$.
Then $$r \le d_P(e,y) <r+m, \,\, r \le  d_P(e,x) < r+m,$$ so that
$d_P(x,y) < m$, and moreover the segment between $x,y$ on $\gamma$
is contained in $A$, therefore $d_A(x,y) < m$. 
It means $x$ is in the open $m$-neighborhood of $p \cup s$ with respect
to $d_A$, which contradicts the way we chose $x$.

{\it Case 2: Planar graphs}.
The argument is almost same as the case 1,
so we will be brief. We also keep the notations. 
Suppose $A$ contains an embedded $m$-fat theta-curve $\Theta(p,q,s)$. $\Theta$ contains only finitely many 
edges, so that $\phi|_{\Theta}$ is continuous.
We proceed as before, and take a geodesic $\gamma$ in $P$.
Again, it contains only finitely many edges, so that 
$\phi|_\gamma$ is continuous and gives a path $\phi(\gamma)$ in $\Bbb R^2$.
So, $\gamma $ must intersect 
$p \cup s$.
The rest is same.
\qed

We will show a few more lemmas.
Although we keep the planar assumption, we only 
use the conclusion of Lemma \ref{annulus.no.theta}, ie,
no embedded, fat theta curves in annuli.

\begin{Lem}\label{lemma.0}
Suppose $(P,d_P)$ is a planar set or a planar graph. 
Given $r, m>0$, let $A$ be a connected component of $A(r,r+5m)$,
and $d_A$ its path metric. Then for any $L>0$ there is a constant $D(L)$, which depends
only on $L$ and $m$, such that  $(A,d_A)$ has a cover by $D(L)$-bounded sets whose 
$L$-multiplicity is at most 2.

Moreover, we can take $D(L)=10^5 \max\{L,5m\}.$ 
\end{Lem}

\proof
 Apply Lemma \ref{annulus.no.theta} to $A$, then 
 $(A,d_A)$ has no
embedded, $5m$-fat theta curve.
Namely, $(A,d_A)$ is a $5m$-coarse cactus. 
 Then, Theorem
\ref{quasi.cactus} implies that a desired constant $D(L)$ exists,
which depends only on $L,m$. The moreover part
is also from the theorem. 
\qed

\subsection{Asymptotic dimension of a plane}

Lemma \ref{lemma.0} implies a similar result with respect to the metric $d_P$
for $L=m$ if 
we reduce the width of the annulus:

\begin{Lem}\label{lemma.1}
Suppose $(P,d_P)$ is a planar set or a planar graph. 
Given $r, m>0$, 
let $A_1(r,r+3m)$ be a connected component of $A(r,r+3m)$. 
Then 
there is a cover of $(A_1(r,r+3m),d_P)$,
by $(10^6m)$-bounded sets whose 
$m$-multiplicity is at most $2$.


\end{Lem}
\proof
Let $A_1(r-m,r+4m)$ be the connected component
of $A(r-m,r+4m)$ that contains $A_1(r,r+3m)$.
Apply the lemma \ref{lemma.0} to $A_1(r-m,r+4m)$ with the path metric, 
setting $L=m$, 
and obtain a cover whose $m$-multiplicity is at most $2$
by $(10^6m)$-bounded sets. Restrict the cover to 
$A_1(r,r+3m)$. We argue this is a desired cover. 
First, this cover is $10^6m$-bounded w.r.t. $d_P$. That is clear
since $d_P$ is not larger than the path metric on $A_1(r-m,r+4m)$.

Also, its $m$-multiplicity is 2 w.r.t. $d_P$.
To see it, let $x \in A_1(r,r+3m)$ be a point. Suppose $K$ is a  set
in the cover with $d_P(x,K) \le m$.
Then a path that realizes the distance $d_P(x,K)$
 is contained in $A_1(r-m,r+4m)$, so that the distance 
 between $x$ and $K$ is at most $m$ w.r.t.
 the path metric on $A_1(r-m,r+4m)$.
 But there are at most 2 such $K$ for a given $x$,
 and we are done. 
\qed

Lemma \ref{lemma.1} implies a lemma for the entire annulus, if we reduce the width further,
 which is in general not connected.

\begin{Lem}\label{lemma.2}
Suppose $(P,d_P)$ is a planar set or a planar graph. 
Then, for any $r, m>0$, there is a cover of $(A(r,r+m),d_P)$
by $(10^6m)$-bounded sets whose 
$m$-multiplicity is at most 2.


\end{Lem}
\begin{figure}[htbp]
\hspace*{-3.3cm}     
\begin{center}
                                                      
\includegraphics[scale=0.18]{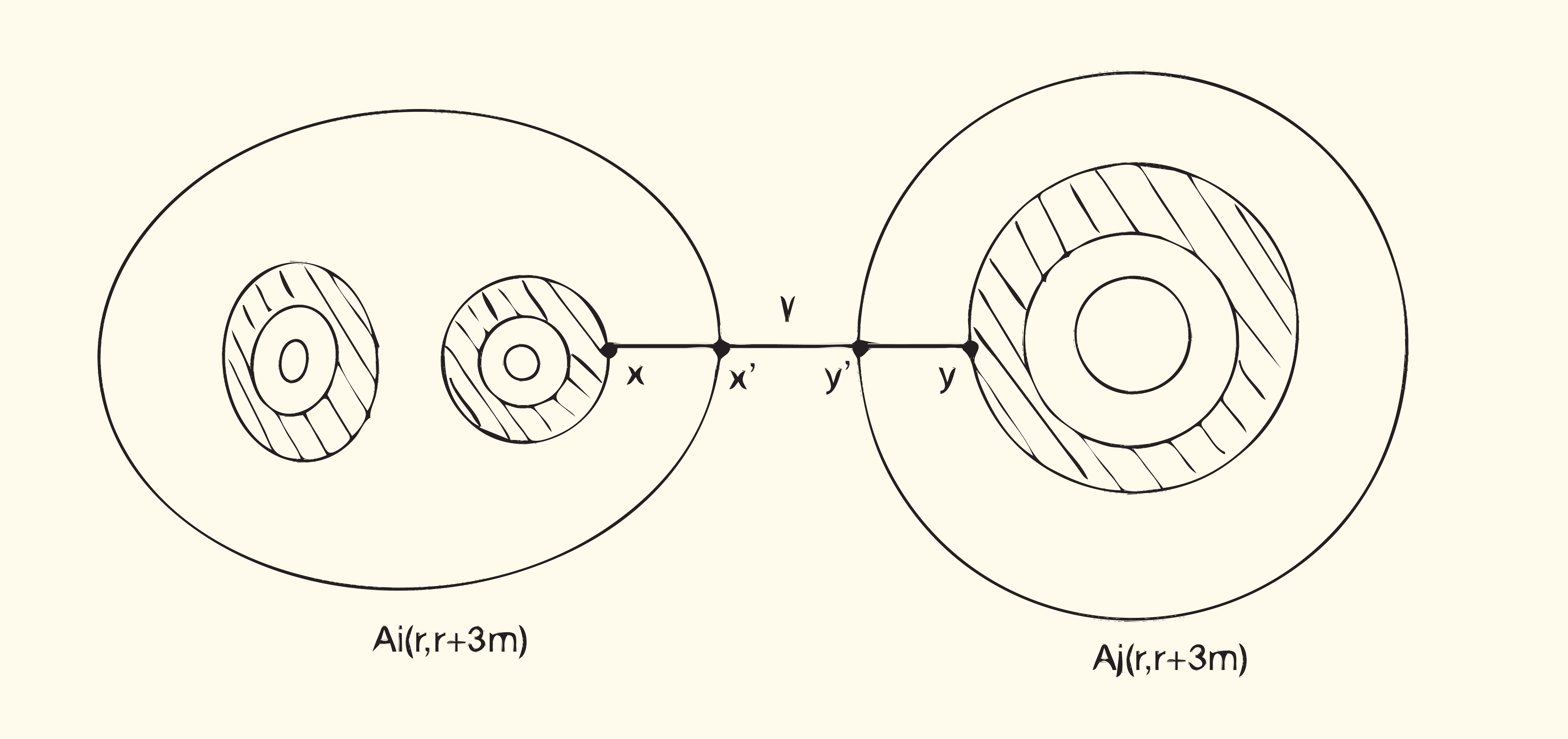}
  \end{center}
  \caption{The shaded area in $A_k(r,r+3m)$ is $A_k(r+m,r+2m)$ for $k=i,j$.
  $[x,x'] \subset A_i(r,r+3m), [y,y']\subset A_j(r,r+3m)$.}
  \label{annuli.cover}
\end{figure}

\proof
We will construct a desired covering for $(A(r+m,r+2m),d_P)$, then 
rename $r+m$ by $r$.
(Strictly speaking, this renaming works only for $r>m$. But if $r\le m$, then 
the diameter of $A(r,r+m)$ is $\le 4m$, so that the conclusion holds.)

The metric in the argument is  $d_P$ unless otherwise said.

Let $A_1(r,r+3m)$ be a connected component of $A(r,r+3m)$.
By  lemma \ref{lemma.1}, we have a covering of $(A_1(r,r+3m), d_P)$
by $(10^6m)$-bounded sets whose $m$-multiplicity is 2. 
Then restrict the covering to the set 
$$A_1(r+m, r+2m)=
A_1(r,r+3m) \cap A(r+m, r+2m).$$

Apply the same argument to all other components,
$A_i(r,r+3m)$, of $A(r,r+3m)$, and 
obtain a covering for 
$$A_i(r+m, r+2m)=
A_i(r,r+3m) \cap A(r+m, r+2m).$$
So far, we obtained a desired covering for each $A_i(r+m, r+2m)$.

Consider the following decomposition, 
$$A(r+m, r+2m) =\sqcup_i A_i(r+m, r+2m).$$
We will obtain a desired covering on the left hand side
by gathering the covering we have for each set on the right hand side. 
We are left to verify that the sets $A_i(r+m,r+2m)$'s are $2m$-separated
from each other w.r.t. $d_P$.

Indeed, let $A_i(r+m,r+2m), A_j(r+m, r+2m)$
be distinct sets.
Then
$$A_i(r+m,r+2m) \subset A_i(r,r+3m), A_j(r+m,r+2m) \subset A_j(r,r+3m),$$
$$A_i(r,r+3m) \cap A_j(r,r+3m) = \emptyset.$$
Now, take a point $x \in A_i(r+m,r+2m)$
and a point $y \in A_j(r+m, r+2m)$.
Join $x,y$ by a geodesic, $\gamma$, in $P$.
See Figure \ref{annuli.cover}.
Let $x' \in \gamma$ be the first point
where $\gamma$ exits $A_i(r,r+3m)$.
Then we have
$d_P(x,x') \ge m.$
Let $y' \in \gamma$ be the last point where
$\gamma$ enters $A_j(r,r+3m)$.
Then
$d_P(y',y) \ge m$.
Since $A_i(r,r+3m)$ and $A_j(r,r+3m)$
are disjoint, 
$$d_P(x,y) > d_P(x,x') + d_P(y',y) = 2m.$$
%
\qed

\subsection{Proof of Theorems \ref{Main}, \ref{inf-graph}
and  Corollary \ref{fin-graph}}
\label{section.proof}
We prove Theorems \ref{Main} and \ref{inf-graph} at one time. 

\proof
By assumption, $(P,d_P)$ is either a planar set (Theorem \ref{Main}) or a planar graph
(Theorem \ref{inf-graph}).
Given $m>0$, define annuli
$$A_n=A(nm, (n+1)m), n\ge 0.$$ Set $D(m)=10^6m$.
By Lemma \ref{lemma.2} each $(A_n,d_P)$ has a covering by $D(m)$-bounded sets
whose $m$-multiplicity is at most 2. 

Gathering all of the coverings  for the annuli, we 
have a covering of $(P,d_P)$ by $D(m)$-bounded sets
whose $\frac{m}{3}$-multiplicity is at most 4
since any ball of radius $\frac{m}{3}$ intersect 
at most two annuli as $A_n$ and $A_{n+2}$ are at least $m$-apart for all $n$
with respect to $d_P$.
We are done by renaming $\frac{m}{3}$ by $m$, and changing 
$D(m)$ to
$D(m)=3(10^6m)$ accordingly. 
\qed

There is nothing more to argue for Corollary \ref{fin-graph} since 
it is only a special case of Theorem \ref{inf-graph}
for finite graphs.

\section{Questions and remarks}
An obvious open question is the following:
\begin{Qu}\label{Q1}
Is the asymptotic dimension of a plane
 at most two for any geodesic metric?
\end{Qu}

{\it Note added in proof}.
J{\o}rgensen-Lang \cite{JL} have answered the question affirmatively 
by now.
An argument goes like this (slightly different from \cite{JL}). 
For a map $f:X \to Y$ between metric spaces, 
Brodskiy-Dydak-Levin-Mitra \cite{BDLM} introduced
the notion of the asymptotic dimension of $f$, $\asdim f$, and 
proved a Hurewicz type theorem, \cite[Theorem 4.11]{BDLM}:
$\asdim X \le \asdim f + \asdim Y$.
Now apply this to the distance function from a base point,
$f: P \to \Bbb R$.
Using Lemma \ref{lemma.2} one argues $\asdim f \le 1$, and 
since $\asdim \Bbb R =1$, it follows $\asdim P \le 2$.
This is only for the asymptotic dimension, and they \cite{JL}
showed the Assouad-Nagata dimension of $P$ is at most 2
by exhibiting  a linear bound for $D(m)$. 
Also, concerning Question \ref{Q1} another proof
of a slightly more general result is given by 
Bonamy-Bousquet-Esperet-Groenland-Pirot-Scott
\cite{BBEGPS}.

It is reasonable to ask whether the asymptotic bound for minor excluded graphs is uniform:

\begin{Qu} \label{Q2}
Given $m\ge3$, is there an $M>0$ such that if $\Gamma$ be a connected graph 
excluding the complete graph $K_m$ as a minor then $\Gamma$
has asymptotic dimension at most $M$? In fact one may ask whether it is possible to take $M=2$.
\end{Qu}

{\it Note added in proof}.
Bonamy et al  \cite{BBEGPS} have answered this by now in the bounded degree case and Liu \cite{Li} in general.

In contrast to Theorem \ref{Main},
\begin{Prop}
$\Bbb R^3$ has a Riemannian metric whose 
asymptotic dimension is infinite.
\end{Prop}

Probably this result is known to experts but we give a proof as 
we did not find it in the literature. Note that any finite graph can be embedded in $\mathbb R^3$ and one sees easily that by changing the metric one can make these embeddings say $(2,2)$ quasi-isometric.
Indeed one may take a small neighborhood of the graph and define a metric so that the distance from an edge to the surface of this neighborhood 
is sufficiently large. Fix $n>3$ and take a unit cubical grid in $\Bbb R^n$, then consider a sequence of finite subgraphs $\Gamma _i$ in the grid of size $i >0$.
 We join $\Gamma _i$ with $\Gamma _{i+1}$ by an edge (for all $i$) and 
we obtain an infinite graph, $\Lambda^n$, whose asymptotic dimension is
equal to $n$. This graph also embeds in $\Bbb R^3$ and 
one can arrange a Riemannian metric on $\Bbb R^3$
such that the embedding is $(2,2)$ quasi-isometric.
For this metric the asymptotic dimension of $\Bbb R^3$
is at least $n$.
Finally we can embed the disjoint union of $\Lambda^n,
n>3$ in $\Bbb R^3$ and arrange a Riemannian metric 
on $\Bbb R^3$ such that the embedding is $(2,2)$ quasi-isometric. Now the asymptotic dimension of $\Bbb R^3$ is infinite for this metric.

\end{document}